\documentclass[a4paper]{article}
\usepackage{latexsym}
\usepackage{amssymb}
\usepackage{amsthm}
\usepackage{amsmath,amsfonts}

\def\titlerunning#1{\gdef\titrun{#1}}
\makeatletter
\def\author#1{\gdef\autrun{\def\and{\unskip, }#1}\gdef\@author{#1}}
\def\address#1{{\def\and{\\\hspace*{18pt}}\renewcommand{\thefootnote}{}%
\footnote {#1}}%
\markboth{\autrun}{\titrun}}
\makeatother
\def\email#1{\hspace*{4pt}{\em e-mail}: #1}

\newcommand{\Ff}{{\mathbb F}}
\newcommand{\adj}{\,\sim\,}
\newcommand{\lspan}[1]{{\langle{#1}\rangle}}
\newcommand{\gcolon}{{\,:\,}}

\titlerunning{}

\title{A construction of directed strongly regular graphs
with parameters (63,11,8,1,2)}
\author{Andries E. Brouwer, Dean Crnkovi\' c and Andrea \v Svob}
\date{2024-04-15}

\begin{document}
\maketitle

\address{D. Crnkovi\' c, A. \v Svob: Faculty of Mathematics, University of Rijeka, Radmile Matej\v ci\'c 2, 51000 Rijeka, Croatia;
\email{aeb@cwi.nl, \{deanc,asvob\}@math.uniri.hr}
}

\begin{abstract}
In this paper, we prove the existence of directed strongly regular graphs with parameters $(63,11,8,1,2)$. We construct a pair of nonisomorphic dsrg(63,11,8,1,2),
where one is obtained from the other by reversing all arrows. Both directed strongly regular graphs have $L_2(8):3$ as the full automorphism group.
\end{abstract}

\bigskip



\section{Introduction}
Given a finite directed graph $\Delta$ with vertex set $X$ and without loops,
let its adjacency matrix be the matrix $A$ with rows and columns indexed
by $X$ such that $A_{xy} = 1$ when $xy$ is an edge, and $A_{xy} = 0$
otherwise.

Let $I$ denote the identity matrix of suitable size,
and $J$ the all-1 matrix.

A directed strongly regular graph (dsrg) with parameters
$v,k,t,\lambda,\mu$ is a directed graph on $v$ vertices,
such that each vertex has indegree and outdegree $k$,
where the adjacency matrix $A$ satisfies the equation
$A^2 = tI + \lambda A + \mu (J-I-A)$.
Equivalently, such that for any two vertices $x,y$ the number of
directed paths $x \to z \to y$ is $t$, $\lambda$ or $\mu$
when $x=y$, or $xy$ is an edge, or $x$ is not an edge, respectively.

This concept was defined by Duval \cite{Duval88} as a directed
generalization of that of strongly regular graph.
See also \cite{BvM}, \S8.22.
Inspection of the literature and the table in \cite{BH} suggests
that no dsrg(63,11,8,1,2) was known. In this note, we construct two nonisomorphic directed strongly regular graphs with parameters (63,11,8,1,2), that are the first known directed strongly regular graphs with this parameters.

\section{A bundle of conics}
For $q = 2^e$, let $F = \Ff_q$ be the finite field of order $q$.
Since squaring is an automorphism of $F$, each element $a$ has
a well-defined square root $a^{1/2}$.

Consider the projective plane $PG(2,q)$ provided with a nondegenerate
qua\-dra\-tic form $Q$ defining a conic $C$ with nucleus $N$.
Let $(,)$ be the bilinear form given by $Q(x+y)=Q(x)+Q(y)+(x,y)$.
Let $X$ be the set of $q^2-1$ nonsingular points other than the nucleus.
Since $Q(ax) = a^2 Q(x)$ for $a \in F$, $x \in X$,
the point $\lspan{x}$ has a unique representation
$\bar{x} = Q(x)^{-1/2} x$.

Let $\Gamma$ be the (undirected) graph with vertex set $X$, where
$\lspan{x} \adj \lspan{y}$ when $(\bar{x},\bar{y}) = 1$. Equivalently, when
$(x,y)^2 = Q(x)Q(y)$. Then $\Gamma$ is a distance-regular
graph with intersection array $\{q,q-2,1;\,1,q-2,q\}$,
an antipodal $(q-1)$-cover of the complete graph $K_{q+1}$,
with $q^2-1$ vertices and $\lambda=\mu=1$.
Two vertices $\lspan{x},\lspan{y}$ are antipodal when $(\bar{x},\bar{y}) = 0$.
See also \cite{BCN}, Proposition 12.5.3.

The automorphism group $G$ of $\Gamma$ is the semilinear group
$P\Sigma L(2,q)$, acting transitively and edge-transitively on $\Gamma$.
The orbits of the stabilizer in $L_2(q)$ of the vertex $p$
(which is elementary abelian of order $q$) are the $q-1$ singletons
in the antipodal class of $p$ together with the $q-1$ sets of size $q$
that are the conics with equations $Q(p)Q(x)+c(x,p)^2 = 0$ for
$c \in \Ff_q \setminus \{0\}$ minus the point $C \cap p^\perp$.

\medskip
$\Gamma$:
{\scriptsize

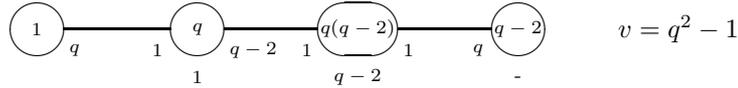
$$\begin{picture}(240,15)(0,-11)
\put(0,0){\circle{20}}\put(0,0){\makebox(0,0){1}}
\put(10,0){\line(1,0){40}}
\put(12,-8){\makebox(0,0)[l]{$q$}}
\put(47,-8){\makebox(0,0)[r]{1}}
\put(60,0){\circle{20}}\put(60,0){\makebox(0,0){$q$}}
\put(60,-18){\makebox(0,0){1}}
\put(70,0){\line(1,0){35}}
\put(72,-8){\makebox(0,0)[l]{$q-2$}}
\put(103,-8){\makebox(0,0)[r]{1}}
\put(120,0){\oval(30,20)}\put(120,0){\makebox(0,0){$q(q-2)$}}
\put(120,-18){\makebox(0,0){$q-2$}}
\put(135,0){\line(1,0){35}}
\put(137,-8){\makebox(0,0)[l]{1}}
\put(167,-8){\makebox(0,0)[r]{$q$}}
\put(180,0){\circle{20}}\put(180,0){\makebox(0,0){$q-2$}}
\put(180,-18){\makebox(0,0){-}}
\put(240,0){\makebox(0,0){\normalsize $v=q^2-1$}}
\end{picture}$$
}

\section{A directed strongly regular graph}
The stabilizer in $L_2(q)$ of one point in an antipodal class
fixes that class pointwise. The stabilizer in $G$ of one point
in an antipodal class has orbits of lengths dividing $e$ (on that class).

In the special case $q = 8$, with $G = L_2(8) \gcolon 3$,
the six antipodes of a point $p$ fall into two orbits of size 3.
Adding to $\Gamma$ (viewed as a directed graph by viewing an
undirected edge $xy$ as a pair of directed edges $xy$ and $yx$)
arrows from $p$ to the three antipodes in one orbit,
and adding all images of these arrows under $G$,
produces a directed strongly regular graph with parameters $(63,11,8,1,2)$.
Let $\Delta_1$ and $\Delta_2$ be the two dsrgs obtained for the two choices
of orbit, with adjacency matrices $A_1$ and $A_2$, respectively.
Then $\Delta_2$ is not isomorphic to $\Delta_1$, but $\Delta_2$ is
isomorphic to the dsrg with adjacency matrix $A_1^\top$.

\bigskip
$\Delta$:
{\scriptsize

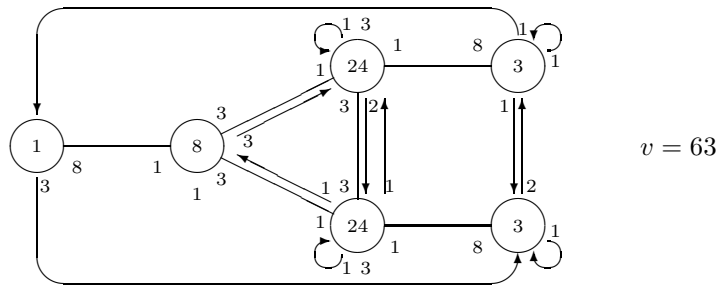
$$\begin{picture}(240,70)(0,-30)
\put(0,0){\circle{20}}\put(0,0){\makebox(0,0){1}}
\put(60,0){\circle{20}}\put(60,0){\makebox(0,0){8}}
\put(120,30){\circle{20}}\put(120,30){\makebox(0,0){24}}
\put(120,-30){\circle{20}}\put(120,-30){\makebox(0,0){24}}
\put(180,30){\circle{20}}\put(180,30){\makebox(0,0){3}}
\put(180,-30){\circle{20}}\put(180,-30){\makebox(0,0){3}}
\put(10,0){\line(1,0){40}}
\put(13,-8){\makebox(0,0)[l]{8}}
\put(47,-8){\makebox(0,0)[r]{1}}
\put(60,-18){\makebox(0,0){1}}
\put(68.944,4.472){\line(2,1){42.111}}
\put(68.944,12.472){\makebox(0,0)[c]{3}}
\put(75,4){\vector(2,1){35}}
\put(79,2){\makebox(0,0)[c]{3}}
\put(108.056,28.528){\makebox(0,0)[r]{1}}
\put(68.944,-4.472){\line(2,-1){42.111}}
\put(68.944,-12.472){\makebox(0,0)[c]{3}}
\put(110,-21){\vector(-2,1){35}}
\put(108,-16){\makebox(0,0)[c]{1}}
\put(108.056,-28.528){\makebox(0,0)[r]{1}}
\put(109,41){\oval(10,10)[t]\oval(10,10)[lb]}
\put(110,35.8){\vector(1,0){0}}
\put(114,46){\makebox(0,0)[l]{1}}
\put(109,-41){\oval(10,10)[b]\oval(10,10)[lt]}
\put(110,-35.8){\vector(1,0){0}}
\put(114,-46){\makebox(0,0)[l]{1}}
\put(117,-15){\makebox(0,0)[r]{3}}
\put(117,15){\makebox(0,0)[r]{3}}
\put(123,46){\makebox(0,0){3}}
\put(123,-46){\makebox(0,0){3}}
\put(120,-20){\line(0,1){40}}
\put(123,18){\vector(0,-1){36}}
\put(130,-18){\vector(0,1){36}}
\put(124,15){\makebox(0,0)[l]{2}}
\put(130,-15){\makebox(0,0)[l]{1}}
\put(130,30){\line(1,0){40}}
\put(133,38){\makebox(0,0)[l]{1}}
\put(167,38){\makebox(0,0)[r]{8}}
\put(130,-30){\line(1,0){40}}
\put(132,-38){\makebox(0,0)[l]{1}}
\put(167,-38){\makebox(0,0)[r]{8}}
\put(178.5,18){\vector(0,-1){36}}
\put(181.5,-18){\vector(0,1){36}}
\put(177,15){\makebox(0,0)[r]{1}}
\put(183,-15){\makebox(0,0)[l]{2}}
\put(191,41){\oval(10,10)[r]\oval(10,10)[lt]}
\put(185.9,39.5){\vector(0,-1){0}}
\put(192,32){\makebox(0,0)[l]{1}}
\put(191,-41){\oval(10,10)[r]\oval(10,10)[lb]}
\put(185.9,-39.5){\vector(0,1){0}}
\put(192,-32){\makebox(0,0)[l]{1}}
\put(90,42){\oval(180,20)[t]}
\put(0,42){\vector(0,-1){30}}
\put(180,44){\makebox(0,0)[l]{1}}
\put(90,-42){\oval(180,20)[b]}
\put(0,-12){\line(0,-1){30}}
\put(1,-15){\makebox(0,0)[l]{3}}
\put(180,-40){\vector(0,1){0}}
\put(240,0){\makebox(0,0){\normalsize $v=63$}}
\end{picture}$$
}

\section{Computational construction}

The directed strongly regular graphs $\Delta_1$ and $\Delta_2$ can be obtained by using \cite[Theorem 3]{cms}.
That construction gives us all simple designs on which a group $G$ acts transitively on the points and blocks, {\it i.e.} if $G$ acts 
transitively on the points and blocks of a simple $1$-design ${\mathcal{D}}$, then ${\mathcal{D}}$ can be obtained by the given method. 

Note that the construction from \cite{cms} gives us $1$-designs, and the incidence matrices of some of these 1-designs may be the adjacency matrices of directed strongly regular graphs.
Since that construction gives all designs having $G$ as an automorphism group acting transitively on points and blocks, it gives us also all directed strongly regular graphs admitting a transitive action of the set of vertices. 

The linear group $L_2(8)$ is the simple group of order 504 and up to conjugation it has exactly one subgroup of order 8, which is isomorphic to the elementary abelian group $E_8$. Using the method given in \cite[Theorem 3]{cms}, by taking $G=L_2(8)$ and the stabilizer of a vertex $G_{\alpha}=E_8$, we show that up to isomorphism, there are exactly two directed strongly regular graphs with parameters $(63,11,8,1,2)$ on which the linear group $L_2(8)$ acts transitively, namely $\Delta_1$ and $\Delta_2$. These directed strongly regular graphs have $L_2(8):3$ as the full automorphism group.

The computations are made by using programmes written for Magma \cite{magma}. The adjacency matrices of the two directed strongly regular graphs are available online at {\tt\verb# http://www.math.uniri.hr/~asvob/DSRG_63_11_8_1_2.txt#}


%

\section*{Acknowledgement}
D. Crnkovi\' c and A. \v Svob were supported by {\rm C}roatian Science Foundation under the projects 4571 and 5713.


\begin{thebibliography}{99}

\bibitem{magma}
W.~Bosma, J.~Cannon, Handbook of Magma Functions, Department of Mathematics, University of Sydney, 1994.
{\tt\verb#http://magma.maths.usyd.edu.au/magma#}.

\bibitem{BCN}
A. E. Brouwer, A. M. Cohen \& A. Neumaier,
{\it Distance-regular graphs},
Springer, 1989.

\bibitem{BH}
A. E. Brouwer \& S. A. Hobart,
{\it Parameters of directed strongly regular graphs},
{\tt\verb#http://homepages.cwi.nl/~aeb/math/dsrg/dsrg.html#}.

\bibitem{BvM}
A. E. Brouwer \& H. Van Maldeghem,
{\it Strongly regular graphs},
Cambridge Univ. Press, 2022.


\bibitem{cms}
D.~Crnkovi\'{c}, V.~Mikuli\'{c}~Crnkovi\'{c}, A.~\v{S}vob, 
On some transitive combinatorial structures constructed from the unitary group $U(3,3)$, 
J. Statist. Plann. Inference 144 (2014) 19--40.

\bibitem{Duval88}
A.~M.~Duval,
{\it A directed graph version of strongly regular graphs},
J. Combin. Th. (A) {\bf 47} (1988) 71--100.

\end{thebibliography}
\end{document}